\input amstex
\documentstyle{amsppt}
\nologo
\define\BaseLineSkip{\baselineskip=13.3pt}  
\NoBlackBoxes
\rightheadtext\nofrills{\eightpoint Generalized Bernoulli-Hurwitz Numbers}
\leftheadtext\nofrills{\eightpoint Yoshihiro \^Onishi}
\hcorrection{13mm} 
\vcorrection{-5mm} 
\define\inbox#1{$\boxed{\text{#1}}$}
\def\fp{\flushpar}

\define\underbarl#1{\lower 1.4pt \hbox{\underbar{\raise 1.4pt \hbox{#1}}}}

\define\tp#1{\negthinspace\ ^t#1}

\define\lr#1{^{\sssize\left(#1\right)}}
\define\br#1{^{\sssize\left[#1\right]}}
\define\do#1{^{\sssize\left<#1\right>}}

\define\C#1#2{{{#1}\choose{#2}}}

\define\dg#1{(d^{\circ}\geqq#1)}
\define\Dg#1#2{(d^{\circ}(#1)\geqq#2)}

\define\Dx#1{\left(\frac{d}{dx}\right)^{#1}}




\font\tenptmbit=cmmib10 
\define\bk#1{\text{\tenptmbit{#1}}}

\font\twelveptbf=cmbx12

\def\oVrule{\vrule width .5pt}
\def\oHrule{\hrule height .5pt}
\def\lbox#1#2#3{\kern 0pt %
\dimen1=#1pt \dimen2=#2pt%
\advance\dimen1 by -1.0pt%
\dimen3=\dimen1%
\advance\dimen3 by -16pt%
\advance\dimen2 by -1.0pt%
\dimen4=\dimen2%
\advance\dimen4 by -20pt%
\vbox to #1pt
{
\hsize #2pt
\oHrule
\hbox to #2pt
{
\vsize \dimen1
\oVrule
\vbox to \dimen1
{
\hsize \dimen2
\vskip 8pt
\hbox to \dimen2
{
\vsize \dimen3
\hskip 10pt
%
%
\vbox to \dimen3{\hsize \dimen4\fp #3\vfil}%
\hfil\hskip 10pt
}%
\vskip 8pt
}%
\oVrule%
}%
\oHrule%
}\kern 0pt%
}%
\def\cbox#1#2#3{\kern 0pt %
\dimen1=#1pt \dimen2=#2pt%
\advance\dimen1 by -1.0pt%
\dimen3=\dimen1%
\advance\dimen3 by -3pt%
\advance\dimen2 by -1.0pt%
\dimen4=\dimen2%
\advance\dimen4 by -3pt%
\vbox to #1pt%
{%
\hsize #2pt
\oHrule%
\hbox to #2pt%
{%
\vsize \dimen1%
\oVrule%
\vbox to \dimen1%
{%
\hsize \dimen2%
\vskip 1.5pt
\hbox to \dimen2%
{%
\vsize \dimen3%
\hskip 1.5pt\hfil
\vbox to \dimen3{\hsize \dimen4\vfil\hbox{#3}\vfil}%
\hfil\hskip 1.5pt
}%
\vskip 1.5pt
}%
\oVrule%
}%
\oHrule%
}\kern 0pt
}%
\def\cboxit#1#2#3{$\hbox{\lower 2.5pt \hbox{\cbox{#1}{#2}{#3}}}$}
\topmatter
\title \nofrills{\twelveptbf  
Theory of Generalized Bernoulli-Hurwitz Numbers \\
         for the Algebraic Functions of Cyclotomic Type}
\endtitle
\author {Yoshihiro \^Onishi} 
\endauthor
\address 3-18-34, Ueda, Morioka, 020-8550, Japan
\endaddress
\email onishi{\@}iwate-u.ac.jp
\endemail
\endtopmatter
\document

\BaseLineSkip

\subheading{1. Introduction}
\fp
Let  $C$  be the projective curve of genus  $g$  defined by 
  $$
  y^2=x^{2g+1}-1 \ \ \ (\text{or} \ \ \ y^2=x^{2g+1}-x).
  \tag 1.1
  $$
The unique point of  $C$  at infinity is denoted by  $\infty$.  
Let consider the integral
  $$
  u=\int_{\infty}^x\frac{x^{g-1}dx}{2y}.
  \tag 1.2
  $$
This is an integral of a differential of first kind on  $C$  
which does not vanish at $\infty$.  
The integral converges everywhere.  
If  $g=0$  the inverse funcion of (1.2) is  $-1/\sin^2(u)$.  
As is well-known, if  $g=1$  the inverse function of (1.2)  is 
just the Weierstrass function  $\wp(u)$  with  $\wp'(u)^2=4\wp(u)^3-4$  
(or  $\wp'(u)^2=4\wp(u)^3-4$).  

The Bernoulli numbers  $\{B_{2n}\}$  are the coefficients of the Laurent expansion 
of $-1/\sin^2(u)$  at $u=0$:
  $$
  \frac{-1}{\sin^2(u)}=\frac{1}{u^2}
  -\sum_{n=1}^{\infty}(-1)^n\frac{2^{2n}B_{2n}}{2n}\frac{u^{2n-2}}{(2n-2)!}.
  \tag 1.3
  $$
The Hurwitz numbers  $\{H_{4n}\}$  are the coefficients of the expansion
  $$
  \wp(u)=\frac1{u^2}
  +\sum_{n=1}^{\infty}\frac{2^{4n}H_{4n}}{4n}\frac{u^{4n-2}}{(4n-2)!}.
  \tag 1.4
  $$
These kinds of numbers are quite important, especially, in number theory.  
Because the inverse function on the neighborhood of  $u=0$   
(1.2) can not extend globally with respect to  $u$ 
if $g>1$, most of mathematician never considered (1.2) for $g>1$,  
usually based on the classical idea of Jacobi, 
and worked in several variable functions.  

Surprisingly, for the case  $g>1$, 
the Laurent coefficients of the inverse function of (1.2) have properties 
which properly resemble von Staudt-Clausen's theorem and Kummer's congruence for Bernoulli numbers
(\cite{{\bf C}}, \cite{{\bf vS}},  \cite{{\bf K}}) 
and such theorems for Hurwitz numbers (\cite{{\bf H2}}, \cite{{\bf L}}).  
To explain these facts is the aim of this paper.  
Detailed and extended exposition should be refered to \cite{{\bf \^O}}.

\vskip 5pt

\subheading{2. Main results}
\fp
Here we describe only for the curve   $C$  defined  by
  $$
  y^2=x^5-1.
  $$
We consider the integrals 
  $$
  u_1=\int_{\infty}^{(x,y)}\frac{dx}{2y}, \ \ \ 
  u_2=\int_{\infty}^{(x,y)}\frac{xdx}{2y}
  \tag 2.1
  $$
of the elements of a natural base of the differentials of first kind.  
Since these integrals are converges, 
there exsists a inverse function  $(x(\bk{u}), y(\bk{u}))$  of  $(x,y) \mapsto \bk{u}=(u_1, u_2)$  on 
the range of all values  $\bk{u}$  of (2.1). 
On a neighborhood of $\bk{u}=(0,0)$, 
the functions $x(\bk{u})$ and  $y(\bk{u})$  are 
functions of the second variable  $u_2$  only.  
Hence we have the following differential equation:
  $$
  \frac{du_2}{dx}=\frac{x}{2y}.  
  \tag 2.2
  $$
If we denote  $\frac{dx}{du_2}=x'(\bk{u})$, then
  $$
  x'(u)=\frac{2y}{x}. 
  \tag 2.3
  $$
After squaring the two sides and substituting $y(u)^2=x(u)^5-1$, 
by removing the denominator, we have
\vskip 3pt
\fp
\lbox{43}{392}{\kern 0pt
  $$
  x(\bk{u})^2x'(\bk{u})^2=4x(\bk{u})^5-4
  \tag 2.4
  $$
\vskip 3pt
}\fp
\vskip 3pt
\fp
This (2.4) is just a good analogy of  
  $$
  \wp''(u)=6\wp(u)^2\qquad (\text{or} \ \  \wp''(u)=6\wp(u)^2-2), 
  $$
obtained by  $\wp'(u)^2=4\wp(u)^3-1$ (or $\wp'(u)^2=4\wp(u)^3-4\wp(u)$).   
Indeed, if we define the numbers  $C_{10n}$  and  $D_{10n}$  by 
  $$
  \aligned
  x(\bk{u})&=\frac{1}{{u_2}^2}+\sum_{n=1}^{\infty}
  \frac{C_{10n}}{10n}\frac{{u_2}^{10n-2}}{(10n)!}, \\
  y(\bk{u})&=\frac{-1}{{u_2}^5}+\sum_{n=1}^{\infty}
  \frac{D_{10n}}{10n}\frac{{u_2}^{10n-5}}{(10n)!}, 
  \endaligned
  \tag 2.5
  $$
then we have two Theorems 2.7 and 2.8 below.  
Here, using the property that
  $$
  x(-\zeta u_1, -\zeta^2 u_2)=\zeta x(u_1, u_2), \ \ \ 
  y(-\zeta u_1, -\zeta^2 u_2)=-y(u_1, u_2), 
  \tag 2.6
  $$
we know that only the terms in (2.5) appear. 
The first Theorem is
\vskip 5pt
\fp
\lbox{187}{392}{\kern 0pt \sl
{\bf Theorem 2.7}\ \ \ For each of  $C_{10n}$  and  $D_{10n}$, 
there exist integers   $G_{10n}$  and   $H_{10n}$  such that 
  $$
  \align
  C_{10n} 
  &=\sum_{\Sb p\equiv 1\mod{5} \\ p-1|10n \endSb}
   \negthinspace\negthinspace\negthinspace
   \dfrac{{\ \ A_p}^{10n/(p-1)}}{p}
   \ \ +\ G_{10n}, \\
  D_{10n} 
  &=\sum_{\Sb p\equiv 1\mod{5} \\ p-1|10n \endSb}
   \negthinspace\negthinspace\negthinspace
   \dfrac{(4!^{-1}\ \text{\rm mod}\ p)\ {A_p}^{10n/(p-1)}}{p}
   \ \ +\ H_{10n};
  \endalign
  $$
where 
  $$
  A_p
  =(-1)^{(p-1)/10}
   \dbinom{(p-1)/2}{(p-1)/10}. 
  $$  
}
\vskip 3pt
\fp
This is a generalization of von Staudt-Clausen's theorem.  
The second theorem is
\vskip 5pt
\fp
\lbox{160}{392}{\kern 0pt \sl
{\bf Theorem 2.8}\ \ \ For any prime  $p\equiv 1$ modulo  $5$,  
and positive integers  $a$  and  $n$  such that   $10n-2\geqq a$,  
if  $(p-1)\not\vert\ 10n$, then 
  $$
  \align
  &\sum_{r=0}^a(-1)^r\dbinom{a}{r}{A_p}^{a-r}
  \frac{C_{10n+r(p-1)}}{10n+r(p-1)}\equiv 0 \mod {p^a}, \\
  &\sum_{r=0}^a(-1)^r\dbinom{a}{r}{A_p}^{a-r}
  \frac{D_{10n+r(p-1)}}{(10n+r(p-1)}\equiv 0 \mod {p^a};
  \endalign
  $$
where 
  $$
  A_p=(-1)^{(p-1)/10}
  \dbinom{(p-1)/2}{(p-1)/10}. 
  $$
}\fp
\vskip 5pt
\fp
This is a natural generalization of Kummer's original congruence in \cite{{\bf K}} 
and of such a congruence for Hurwitz numbers (\cite{{\bf L}}, p.190, (20)).  
To prove Theorem 2.8, we need
\vskip 5pt
\fp
\lbox{41}{392}{\kern 0pt \sl
{\bf Theorem 2.9}\ \ \  
For any prime number  $p\equiv 1$ mod $5$   and any positive integer  $n$, 
if $p-1 {\not|}\ 10n$  then  $C_{10n}/10n$  and  $D_{10n}/10n$  belong to  $\bold Z_{(p)}$. 
}
\vskip 3pt
\fp

\vskip 5pt

\subheading{3. About Proof of Theorems} 
\fp
We could not use any addition or multiplication formulae in 
the classical theory of Abelian or Jacobin varieties.  
To prove Therem 2.7 we use a technique based on a method of L. Carlitz \cite{{\bf Car1}}, 
and to prove Theorem 2.9 we use a new technique that is an improved method from that of Carlitz.  
Theorem 2.8 is rather easily shown by Theorem 2.8.  
By using his method, Carlitz succeeded to prove Hurwitz's theorem out side the prime $2$ 
in his paper \cite{{\bf Car1}}.  
Furthermore, although he was trying to find a generalization of von Staudt-Clausen's type theorem, 
Kummer's type congruence for hyperelliptic functions, 
he could not succeed.  
Our results is not only for hyperelliptic curve but also 
for any algebraic curves of type
  $$
  y^a=x^b-1,\ \ \ \text{or}\ \ \ y^2=x^b-x
  $$
with $\gcd(a,b)=1$.  

The detailed proof and a lot of numerical examples are given in \cite{{\bf \^O}}.  


\enddocument
\bye